\documentclass[12pt]{article}
\usepackage{amsmath,amssymb,amsthm,amsfonts,latexsym,hyperref,color}

\usepackage[hmargin=.8in,vmargin=1.5in]{geometry}

\newtheorem{thm}{Theorem}[section]
\newtheorem{prop}[thm]{Proposition}
\newtheorem{cor}[thm]{Corollary}
\newtheorem{lem}[thm]{Lemma}
\newtheorem{conj}[thm]{Conjecture}
\newtheorem{exa}[thm]{Example}

\newcommand{\prf}{\noindent\pf}
\newcommand{\CS}{c^{S}}

\newcommand{\ben}{\begin{enumerate}}
\newcommand{\een}{\end{enumerate}}
\newcommand{\ble}{\begin{lem}}
\newcommand{\ele}{\end{lem}}
\newcommand{\bth}{\begin{thm}}
\renewcommand{\eth}{\end{thm}}
\newcommand{\bpr}{\begin{prop}}
\newcommand{\epr}{\end{prop}}
\newcommand{\bco}{\begin{cor}}
\newcommand{\eco}{\end{cor}}
\newcommand{\bcon}{\begin{conj}}
\newcommand{\econ}{\end{conj}}
\newcommand{\bde}{\begin{defn}}
\newcommand{\ede}{\end{defn}}
\newcommand{\bex}{\begin{exa}}
\newcommand{\eex}{\end{exa}}
\newcommand{\barr}{\begin{array}}
\newcommand{\earr}{\end{array}}
\newcommand{\btab}{\begin{tabular}}
\newcommand{\etab}{\end{tabular}}
\newcommand{\beq}{\begin{equation}}
\newcommand{\eeq}{\end{equation}}
\newcommand{\bea}{\begin{eqnarray*}}
\newcommand{\eea}{\end{eqnarray*}}
\newcommand{\bal}{\begin{align*}}
\newcommand{\bce}{\begin{center}}
\newcommand{\ece}{\end{center}}
\newcommand{\bpi}{\begin{picture}}
\newcommand{\epi}{\end{picture}}
\newcommand{\bpp}{\begin{picture}}
\newcommand{\epp}{\end{picture}}
\newcommand{\bfi}{\begin{figure} \begin{center}}
\newcommand{\efi}{\end{center} \end{figure}}
\newcommand{\bprf}{\begin{proof}}
\newcommand{\eprf}{\end{proof}\medskip}

\newcommand{\bsl}{\begin{slide}{}}
\newcommand{\esl}{\end{slide}}
\newcommand{\bfr}{\begin{frame}}
\newcommand{\efr}{\end{frame}}

\newcommand{\pf}{{\bf Proof.\hspace{5pt}}}

\newcommand{\hqedm}{\hfill \qed \medskip}

\newcommand{\hso}[1]{\hspace{-1pt}}

\newcommand{\emp}{\emptyset}

\newcommand{\sbe}{\subseteq}

\newcommand{\case}[4]{\left\{\barr{ll}#1&\mbox{#2}\\#3&\mbox{#4}\earr\right.}

\newcommand{\flf}[2]{\left\lfloor\frac{#1}{#2}\right\rfloor}
\newcommand{\cef}[2]{\left\lceil\frac{#1}{#2}\right\rceil}

\def\<{\langle}
\def\>{\rangle}

\newcommand{\ree}[1]{(\ref{#1})}

\newcommand{\ka}{\kappa}

\newcommand{\bbN}{{\mathbb N}}

\newcommand{\bbQ}{{\mathbb Q}}

\newcommand{\cB}{{\cal B}}

\newcommand{\fS}{{\mathfrak S}}

\newcommand{\dil}{\displaystyle}

\begin{document}
\pagestyle{plain}

\title{Permutations with given peak set
}
\author{
Sara Billey \thanks{Partially supported by grant DMS-1101017 from
the NSF.} \\[-5pt]
\small Department of Mathematics, University of Washington,\\[-5pt]
\small Seattle, WA 98195-4350, USA, {\tt billey@math.washington.edu}\\
Krzysztof Burdzy \thanks{Partially supported by grant DMS-1206276 from the NSF and by
 grant N N201 397137 from the MNiSW, Poland.
}\\[-5pt]
\small Department of Mathematics, University of Washington,\\[-5pt]
\small Seattle, WA 98195-4350, USA, {\tt burdzy@math.washington.edu}\\
Bruce E. Sagan\\[-5pt]
\small Department of Mathematics, Michigan State University,\\[-5pt]
\small East Lansing, MI 48824-1027, USA, {\tt sagan@math.msu.edu}
}

\date{\today\\[10pt]
	\begin{flushleft}
	                                       \\[5pt]
	\end{flushleft}}

\maketitle

\begin{abstract}
Let $\fS_n$ denote the symmetric group of all permutations
$\pi=a_1\ldots a_n$ of $\{1,\ldots,n\}$.  An index $i$ is a
\emph{peak} of $\pi$ if $a_{i-1}<a_i>a_{i+1}$ and we let $P(\pi)$ be
the set of peaks of $\pi$.  Given any set $S$ of positive integers we
define $P(S;n)$ to be the set $\pi\in\fS_n$ with $P(\pi)=S$.  Our main
result is that for all fixed subsets of positive integers $S$
and all sufficiently large $n$ we have $\# P(S;n)=p(n)2^{n-\#S-1}$ for
some polynomial $p(n)$ depending on $S$.  We explicitly compute $p(n)$
for various $S$ of probabilistic interest, including certain
cases where $S$ depends on $n$.  We also discuss two conjectures, one
about positivity of the coefficients of the expansion of $p(n)$ in a
binomial coefficient basis, and the other about sets $S$ maximizing
$\# P(S;n)$ when $\#S$ is fixed.
\end{abstract}

\section{Introduction}

Let $\bbN$ be the nonnegative integers and, for $n\in\bbN$, let $[n]=\{1,\ldots,n\}$.  Consider the \emph{symmetric group} $\fS_n$ of all permutations $\pi=a_1\ldots a_n$ of $[n]$.  Call an index $i$ a \emph{peak} of $\pi$ if $a_{i-1}<a_i>a_{i+1}$ and define the \emph{peak set} of $\pi$ to be
$$
P(\pi)=\{i\ :\ \mbox{$i$ is a peak of $\pi$}\}.
$$
By way of illustration, if $\pi =a_1\ldots a_7=1453276$ then
$P(\pi)=\{3,6\}$ because of $a_3=5$ and $a_6=7$.  Note that some
authors call $a_i$ a peak rather than $i$.  But our convention is more
consistent with what is used for other permutation statistics such as
the one for descents.  Also note that if $\pi\in\fS_n$ then
$P(\pi)\sbe\{2,\ldots,n-1\}$.  There has been a great deal of research
into peaks of permutations, much of it focusing on permutations with a
given number of peaks.  For example,  
see~\cite{ano:nrp,BHvW,bcmyy:vpp,fv:psl,km:trs,ma:dpe,nym:pas,pet:epp,schocker,ste:epp,str:eap,ws:pen}.

The purpose of the present work is to investigate permutations with  a given peak set.  To this end, define
$$
P(S;n) =\{\pi\in\fS_n\ :\ P(\pi)=S\}.
$$
We will omit the curly brackets around $S$ in this notation.  So, for example,
$$
P(2;4)=\{1324,1423,1432,2314,2413,2431,3412,3421\}.
$$

Our main result will be about the cardinality $\# P(S;n)$ as $n$
varies, where $S$ is a set of constants not depending on $n$.  To
state it, define a set $S=\{i_1<\cdots<i_s\}$ to be
\emph{$n$-admissible} if $\#P(S;n)\neq 0$.  Note that we insist
the elements be listed in increasing order.  The minimum
possible value of $n$ for which $S$ is $n$-admissible is $i_{s}+1$,
and in that case $S$ is $n$-admissible for all $n\ge i_S+1$.  If
we make a statement about an \emph{admissible} set $S$, we
mean that $S$ is $n$-admissible for some $n$ and 
the statement holds
for every $n$ such that $S$ is $n$-admissible.  We can now state our
principal theorem.

\bth \label{main} If $S=\{i_1<\cdots<i_s\}$ is admissible then
$$
\#P(S;n) = p(n)2^{n-\#S-1} 
$$
where $p(n)=p(S;n)$ is a polynomial depending on $S$ such that $p(n)$ is an integer for all integral $n$.  In addition,  
$\deg p(n) = i_s - 1$ (when $S=\emp$ we have $\deg p(n)=0$). 
\eth

If $S$ is not admissible, then $\#P(S;n)=0$ for all positive
integers $n$, so we define the corresponding polynomial $p(S;n)=0$.
Thus, for all sets $S$ of constants not depending on $n$, $p(S;n)$ is
a well defined polynomial which we will call the  \emph{peak
polynomial} for $S$ and $\#P(S;n)=p(S;n)2^{n-\#S -1}$ for all $n>\max
S$ if $S\neq \emptyset$   or for all  $n\geq 1$ if $S=\emptyset$.

Some of the motivation for our work comes from probability theory.  A relationship between permutations and random data has been noticed for quite some time. We refer the reader to the 1937  paper of Kermack and McKendrick~\cite{km:trs} and the references therein. Many probabilistic models are concerned with i.i.d. (independent identically distributed) sequences of data, or their generalization, exchangeable sequences. By definition, any permutation of an exchangeable sequence of data is as likely to be observed as the original sequence. One way to test whether a given sequence of $n$ data points is in fact exchangeable is to analyze the order in which the data are arranged, starting from the highest value to the lowest value. Under the assumption of exchangeability, the order should be a randomly (uniformly) chosen permutation of   $[n]$. Hence, probabilists are interested in probabilities of various events related to uniformly chosen permutations. This is equivalent to evaluating cardinalities of various subsets of $\fS_n$. This article is inspired by, and provides estimates for, a probabilistic project concerned with mass redistribution, to be presented in a forthcoming article~\cite{bbps:me}.  The reader can consult this paper for details, including the specific  i.i.d. sequence which will be used in our model.

The rest of this paper is organized as follows.  The following section
is devoted to a proof of Theorem~\ref{main} and its enumerative
consequences.  The next two sections are devoted to computing the
polynomial $p(n)$ for various sets of interest for the probabilistic
applications. Section~\ref{pc} investigates a conjecture about the
expansion of $p(S;n)$ in a binomial coefficient basis for the space of
polynomials.  The following section states a conjecture about which
$S$ maximizes $\# P(S;n)$ among all $S$ with given cardinality.
Section~\ref{pv} shows how our methods can be applied to the
enumeration of permutations with a fixed set of peaks and valleys.
Finally, we end by using our results to prove a known formula for the
number of permutations with a given number of peaks and suggest
an avenue for future research.

\section{The main enumeration theorem}
\label{met}

We need the following result as a base case for induction.

\bpr
\label{P(emp)}
For $n\ge1$ we have
$$
\#P(\emp;n)=2^{n-1}
$$
\epr
\prf If $\pi\in P(\emp;n)$ then write $\pi=\pi_1 1\pi_2$ where $\pi_1,\pi_2$ are the portions of $\pi$ to the left and right of $1$, respectively.  Now $P(\pi)=\emp$ if and only if $\pi_1$ is decreasing and $\pi_2$ is increasing.  So $\#P(\emp;n)$ is the number of choices of a subset of elements from $[2,n]$ to be in $\pi_1$ since after that choice is made the rest of $\pi$ is determined.  The result follows.\hqedm

We now prove our principal theorem, restating it here for ease of reference.

\bth
\label{P(S):th}
If  $S=\{i_1<\cdots<i_s\}$ is admissible then
$$
\#P(S;n) = p(n)2^{n-\#S-1} 
$$
where $p(n)=p(S;n)$ is a polynomial depending on $S$ such that $p(n)$ is an integer for all integral $n$.  In addition, 
 $\deg p(n) = i_s - 1$ (when $S=\emp$ we have $\deg p(n)=0$). 
\eth

\prf 
We induct on $i=i_1+\cdots+i_s$.  By Proposition~\ref{P(emp)} the
result is true when $i=0$.  Now suppose $i\neq 0$.  For ease of
notation, let $k=i_s-1$ and $S_1=S-\{i_s\}$.  For any fixed
$n> i_{s}$, consider the set $\Pi$ of permutations $\pi=a_1\ldots
a_n\in\fS_n$ such that $P(a_1\ldots a_k)=S_1$ and $P(a_{k+1}\ldots
a_n)=\emp$.  
Since $S$ is $n$-admissible, we know $S_{1}$ is also
$n$-admissible.  Thus, $\#\Pi \neq 0$.

We can construct $\Pi$ by first picking the set of elements to be used for $a_1,\ldots,a_k$ and then arranging this set and its complement to have the prescribed peak sets.  Thus, by induction, the total number of choices is
$$
\#\Pi={n\choose k}P(S_1;k)P(\emp;n-k)={n\choose k} p_1(k) 2^{k-s}\cdot 2^{n-k-1}=
p_1(k){n\choose k} 2^{n-s-1}
$$
for some polynomial $p_1(n)$ with $\deg p_1(n)=i_{s-1}-1<k$.  Also, $p_1(k)$ is an integer which must also be nonzero  since $\#\Pi\neq 0$.

On the other hand, we can count $\Pi$ as follows.  Let $S_2=S_1\cup\{i_s-1\}$.  Note that by the restrictions on $P(a_1\ldots a_k)$ and $P(a_{k+1}\ldots a_n)$ we must have $P(\pi)=S_1$, $P(\pi)=S_2$, or $P(\pi)=S$ for all $\pi\in\Pi$.  And by construction, all such $\pi$ appear which shows that we must have the decomposition $\Pi=P(S_1;n)\cup P(S_2;n)\cup P(S;n)$.  Taking cardinalities and applying induction as well as the previous count for $\#\Pi$ yields
\beq
\label{Pi}
\#P(S;n)
=p_1(k){n\choose k} 2^{n-s-1}-p_1(n)2^{n-s}-p_2(n)2^{n-s-1}
=\left[p_1(k){n\choose k}-2p_1(n)-p_2(n)\right] 2^{n-s-1}
\eeq
where $p_2(n)$ is a polynomial in $n$ of degree $i_s-2<k$ which is integral at integers if $S_2$ is admissible and $0$ otherwise.
To complete the proof, note that $p_1(k){n\choose k}$ is a polynomial in $n$ of degree $k$ while $2p_1(n)$ and $p_2(n)$ have degree smaller than $k$.  So the coefficient of $2^{n-s-1}$ above is a polynomial of degree $k=i_s-1$.    We also have that ${n\choose k}$, $2p_1(n)$,and $p_2(n)$ all have integral values at integers, and we have previously established that $p_1(k)$ is an integer.  So the same is true of the difference above.
\hqedm

Examining the proof of
Theorem~\ref{P(S):th}, one sees that it goes through if $n=i_s$ in the
sense that we will obtain $p(S;i_s)=0$ using the fact that
there are no permutations in $P(S;i_s)$.  (On the other hand, we may
have $p(S;l)\neq0$ for $l<i_s$ and so some bound is still needed.)

Equation~\ree{Pi} immediately yields the following recursive formula
for this polynomial.

\bco
\label{P(S):co}
If $S\neq\emp$ is admissible and $m=\max S$ then  
\beq
\label{P(S;n)}
p(S;n)=p_1(m-1){n\choose m-1}-2p_1(n)-p_2(n)
\eeq
where $S_1=S-\{m\}$, $S_2=S_1\cup\{m-1\}$, and $p_i(n)=p(S_i;n)$ for $i=1,2$.\hqedm
\eco

It is well known that any sequence given by a polynomial of
degree $k$ can be completely determined by any consecutive $k+1$
values by the method of finite differences.  See \cite[Sections 2.6 and
5.3]{GKP} or \cite[Proposition 1.9.2]{sta:ec1}.  If $f(n)$ is a polynomial
function of $n$, then define another polynomial function
of $n$ by $\Delta f(n)= f(n+1)-f(n)$. Similarly, $\Delta^{k}f$ is a
polynomial function of $n$ obtained by applying $\Delta$ successively
$k$ times.  Then expressing $f(n)$ in the basis ${n \choose k}$, we
have
\[
f(n) = \sum_{k=0}^{d} \Delta^{k}f (0) \cdot  {n \choose k}
\]
where $d$ is the degree of $f(n)$.  Thus, Theorem~\ref{P(S):th}  gives
a way to find an explicit formula for $\#P(S;n)$ and $p(S;n)$ for any
admissible set $S$.

For example, if $S=\{2,5 \}$ then $p(S,n)$ has degree $4$.  Thus, we
can find $p(S,n)$ from the sequence $\#P(S,n)/2^{n-3}$ for
$n=6,7,8,9,10$.  Either by hand or computer we find
$\#P(S,6)/2^{3}=10$, $\#P(S,7)/2^{4}=35$, $\#P(S,8)/2^{5}=84$, $\#P(S,9)/2^{6}=168$,
$\#P(S,10)/2^{7}=300$, etc.  Taking successive consecutive differences 5 times gives the difference table 
\[
\begin{matrix}
10 & 35 & 84 & 168 & 300 & 495 & 770 & 1144 & 1638\\
25 & 49 & 84 & 132 & 195 & 275 & 374 & 494\\
24 & 35 & 48 & 63 & 80 & 99 & 120 & \\
11 & 13 & 15 & 17 & 19 & 21\\
2 & 2 & 2 & 2 & 2 & \\
0 & 0 & 0 & 0 & 
\end{matrix}.
\]
Therefore, 
\begin{align}\label{e:2.n-1}
p(S,n)&= 10 {n -6\choose 0} + 25 {n -6\choose 1} + 24 {n -6\choose 2} + 11 {n -6\choose 3} + 2 {n -6\choose 4} \\
&= \frac{1}{12}\,n \left( n-5 \right)  \left( n-2 \right)  \left( n-1 \right).
\end{align}
We used the shifted basis ${n -6\choose k}$ here since the smallest
value of $n$ for which $S$ is $n$-admissible is $n=6$ and we want the
sequences aligned properly.

\bco
\label{c:gf}
If $S$ is a nonempty admissible set and $m=\max S$, then $\#P(S;n)$
has a rational generating function of the form
\[
\sum_{n\geq 1} \#P(S;n) x^{n} = \frac{r(x)}{(1-2x)^{m}}
\]
where $r(x)=r(S;x)$ is the polynomial 
\[
r(x)= (1-2x)^{m} \sum_{n\geq 1} \#P(S;n) x^{n} =
\sum_{k=m+1}^{2m-1} x^k \sum_{j=0}^{k-m-1} (-2)^j{m\choose j}  \#P(S;k-j).
\]
Furthermore, for $n\geq 2m$ the following linear recurrence relation
holds
\[
\sum_{j=0}^{m} (-2)^{j}{m \choose j} \#P(S;n-j) =0.
\]
\eco

\prf
These claims follow directly  by applying the theory of rational generating functions in \cite[Theorem 4.1.1 and
Corollary 4.2.1]{sta:ec1} to Theorem~\ref{P(S):th} and the discussion of the $n=m$ case directly following it.\hqedm

Continuing with the example $S=\{2,5\}$, we get the generating function 
\[
\sum_{n\geq 1} \#P(S;n) x^{n} =
\frac{80 x^6-240 x^7+288 x^8-128 x^9}{(1-2 x)^5} = \frac{16 x^6 (1-x) (5-10 x+8 x^2)}{(1-2 x)^5}
\]
and recurrence relation
\[
\#P(S;n) = 10 \#P(S;n-1) - 40 \#P(S;n-2) + 80 \#P(S;n-3)- 80 \#P(S;n-4) + 32 \#P(S;n-5) 
\]
which holds for $n\geq 10$.

\section{Specific peak sets with constant elements}
\label{sps}

We now derive formulas for $\#P(S;n)$ for various sets $S$ of elements which do not vary with $n$.  These will be useful in proving the results needed for probabilistic applications.  
Peaks represent sites with no mass in the mass redistribution model analyzed in \cite{bbps:me}. In that paper, we will use the results of this section to study the distance between ``empty'' sites.

Before stating the equations, we would like to indicate how they were
originally derived as the proofs below are ones given in hindsight.
Note that Corollary~\ref{P(S):co} expresses $p(S;n)$ in terms of
polynomials for peak sets having smaller maxima than $S$.  So by
iterating this recursion, one can find an expression for $p(S;n)$
whose main contribution is from an alternating sum $\sum_k (-1)^k a_k
{n\choose k}$ for certain coefficients $a_k\ge0$.  By next using the
binomial recursion iteratively, one achieves substantial cancellation.
The simplified formula can then be proved directly using
Corollary~\ref{P(S):co} and these are the results and proofs given
below.

\bigskip

\bth
\label{P(m)}
If $S=\{m\}$ is admissible then
$$
p(S;n)={n-1\choose m-1}-1.
$$
\eth
\prf
We induct on $m$ and use the notation of Corollary~\ref{P(S):co}.  If $m=2$ then $S_1=\emp$ and $S_2=\{1\}$.  By Proposition~\ref{P(emp)} we have $p_1(n)=1$.  Also, $p_2(n)=0$ since $S_2$ is not admissible for any $n$.  Now applying Corollary~\ref{P(S):co} gives
$$
p({2};n)=1\cdot{n\choose 1}-2\cdot 1=n-2={n-1\choose 1}-1.
$$

The induction step is similar, except that now $S_2=\{m-1\}$ which is admissible, but whose polynomial is known by induction.  In particular
$$
p({m};n)={n\choose m-1}-2-\left[{n-1\choose m-2}-1\right]={n-1\choose m-1}-1
$$
as desired.
\hqedm

\bth
If $S=\{2,m\}$ is admissible then
$$
p(S;n)=(m-3){n-2\choose m-1}+(m-2){n-2\choose m-2}-{n-2\choose 1}.
$$
\eth
\prf
 The proof is much like the previous proposition where $S_1=\{2\}$ and $S_2=\{2,m-1\}$.  The details are left to the reader.
\hqedm

From this theorem, we immediately get for $S=\{2,5\}$ that 
$$
p(S;n)=2{n-2\choose 4}+3{n-2\choose 3}-{n-2\choose 1}=\frac{1}{12}\,n \left( n-5 \right)  \left( n-2 \right)  \left( n-1 \right).
$$
Note that this is in agreement with our computations in the previous section.

Using the same technique, one can prove the following result whose demonstration is omitted.

\bth
If $S=\{2,m,m+2\}$ is admissible then
$$
p(S;n)=m(m-3){n\choose m+1}-2(m-3){n-2\choose m-1}-2(m-2){n-2\choose m-2}+2{n-2\choose 1}.
$$
\eth

\section{Peak sets depending on $n$}
\label{ps2}

Consider sets of the form
$$
S=\{i_1<i_2<\dots<i_s<n-j_t<\dots<n-j_2<n-j_1\}
$$
where $i_1,\ldots,i_s,j_1\ldots,j_t$ are constants.  Using exactly the same definition of admissibility as with sets of constants, one can show that Theorem~\ref{P(S):th} continues to hold for such $S$. The specific statement is as follows.

\bth 
Let $S=\{i_1<i_2<\dots<i_s<n-j_t<\dots<n-j_2<n-j_1\}$ be admissible. 
Then
$$
\#P(S;n) = p(n)2^{n-\#S-1} 
$$
where $p(n)=p(S;n)$ is a polynomial depending on $S$ such that $p(n)$ is an integer for all integral $n$.  In addition,  

\hspace{150pt}$
\deg p(n)= \left\{ \barr{ll}
0	&\mbox{if $s=t=0$,}\\     
i_s-1	&\mbox{if $s>0$ and $t=0$,}\\
j_t	&\mbox{if $s=0$ and $t>0$,}\\
i_s+j_t-1 &\mbox{else.\hspace*{150pt}\qed}
\earr
\right.
$
\eth

\medskip

 The demonstration is an induction on $i_1+\dots+i_s+j_1+\dots+j_t$ which is similar to the one given previously and so is omitted.

\bigskip

Next, we will use the results of the previous section to obtain formulas for peak sets $S$ with $2,n-1\in S$
which are useful  probabilistically.  
Of all peak sets that a probabilist might consider, two peaks at the (almost) extreme points of a sequence, namely at $2$ and $n-1$, are of the greatest interest because they represent the distribution of the distance between adjacent empty sites in the mass redistribution model in~\cite{bbps:me}. So a peak set which contains $2$ and $n-1$ represents the joint distribution of several consecutive empty sites.

For sets containing $2$ and $n-1$,
 there is an alternative way to compute $p(S;n)$ which is simpler because it avoids alternating sums.
 If $\pi=a_1\ldots a_n\in P(S;n)$ where $2,n-1\in S$ then $n=a_{i_j}$ for some $i_j\in S$.  (If we do not have the hypothesis on $S$, similar reasoning can be applied, but one needs to worry about the possibility that $a_1=n$ or $a_n=n$.)   
Note also that if we consider the reversal 
$\pi^r=a_n\ldots a_1$ then $P(\pi^r)=n+1-S$ where $n+1-S=\{n+1-i_s,\ldots,n+1-i_1\}$ for $S=\{i_1,\ldots,i_s\}$. 
So if we write $\pi=\pi_L n\pi_R$ we have $\pi_L\in P(S_L;i_j-1)$ and $\pi_R^r\in P(S_R^r;n-i_j)$ where 
$S_L=\{i_1,\ldots,i_{j-1}\}$ and $S_R^r=n+1-\{i_{j+1},\ldots,i_s\}$.
These observations yield the recursion
$$
\#P(S;n)=\sum_{j=1}^s \#P(S_L;i_j-1)\#P(S_R^r;n-i_j) {n-1\choose i_j-1}.
$$
Using Theorem~\ref{P(S):th} and canceling powers of 2 gives
\beq
\label{P(S):eq}
2p(S;n)=\sum_{j=1}^s p(S_L;i_j-1) p(S_R^r;n-i_j) {n-1\choose i_j-1}.
\eeq
Because of the complexity of the formulas, we will often keep the $2$ above on the left-hand side of the equation.

\bth
\label{p(2,n-1)}
If $S=\{2,n-1\}$ is admissible then
$$
p(S;n)=(n-1)(n-4).
$$
\eth
\prf 
In this case, equation~\ree{P(S):eq} has two terms.  In the first $S_L=\emp$ and $S_R^r=\{2\}$, while in the second $S_L=\{2\}$ and $S_R^r=\emp$.  Applying Proposition~\ref{P(emp)} and Theorem~\ref{P(m)}, we obtain
$$
2p(S;n)=1\cdot (n-4)\cdot{n-1\choose 1}+(n-4)\cdot1\cdot{n-1\choose 1}
$$
from which the desired equation follows.
\hqedm

\bth
If $S=\{2,m,n-1\}$ is admissible then
\bea
2 p(S;n)
&=&(m-3)(n-m-2){n-1\choose m-1}+(n-1)\left[(m-3){n-4\choose m-1}+(m-2){n-4\choose m-2}\right.\\[5pt]
&&\hspace*{100pt} \left.+(n-m-1){n-4\choose m-3}+(n-m-2){n-4\choose m-4}-2{n-4\choose1}\right].
\eea
For fixed $n$, the sequence $p(S;n)$ as $m$ varies is symmetric and unimodal and only 
attains its maximum at $m=\flf{n+1}{2}$ and at $m=\cef{n+1}{2}$.
\eth
\prf
The formula for $2p(S;n)$ follows from equation~\ree{P(S):eq} and the results of the previous section similarly to the proof of Theorem~\ref{p(2,n-1)}.  So we leave the details to the reader.

For fixed $n$, let us write $f(m)=2 p(2,m,n-1;n)$ where $4\le m\le n-3$.  The fact that this sequence is symmetric follows directly from the form of $S$.  To prove the rest of the theorem, it suffices to show that the first half of the sequence is strictly increasing.  So consider the difference $f(m+1)-f(m)$ where $m\le (n-1)/2$.  The first term of $f(m)$ contributes
$$
(m-2)(n-m-3){n-1\choose m}-(m-3)(n-m-2){n-1\choose m-1}>0
$$
since, for the given range of $m$, we have $(m-2)(n-m-3)\ge (m-3)(n-m-2)$ by log concavity of the integers and 
${n-1\choose m}>{n-1\choose m-1}$ by unimodality of the binomial coefficients.  Now consider the  the terms with a factor of $n-1$.  Combining terms corresponding to the same binomial coefficient and then using binomial coefficient unimodality gives a contribution to the difference of
$$
\barr{l}
\dil(m-2){n-4\choose m}+2{n-4\choose m-1}+(n-2m){n-4\choose m-2}-2{n-4\choose m-3}-(n-m-2){n-4\choose m-4}\\[20pt]
\dil\hspace*{20pt}>[(m-2)+(n-2m)-(n-m-2)]{m-4\choose m-4}+2\left[{n-4\choose m-1}-{n-4\choose m-3}\right]\\[20pt]
\dil\hspace*{20pt}=2\left[{n-4\choose m-1}-{n-4\choose m-3}\right]\\[20pt]
\dil\hspace*{20pt}>0
\earr
$$
which is what we wished to show.
\hqedm

The demonstration of the next theorem contains no new ideas and so is omitted.

\bth
If $S=\{2,m,m+2,n-1\}$ is admissible then
\bea
2 p(S;n)&=&(m-3)(n-m-4)\left[m{n-1\choose m+1}+(n-m-1){n-1\choose m-1}\right]\\[5pt]
&&\hspace*{50pt} +(n-1)\left[m(m-3){n-2\choose m+1}+(n-m-1)(n-m-4){n-2\choose m-2}\right.\\[5pt]
&&\hspace*{110pt} \left. -2(n-6){n-4\choose m-1} -2(n-6){n-4\choose m-2}+4{n-4\choose 1}\right].
\eea
\eth

\section{A positivity conjecture}
\label{pc}

Given any integer $m$ we have the following basis for the polynomials in $n$
$$
\cB_m=\left\{{n-m\choose k}\ :\ k\ge0\right\}.
$$
Consider a polynomial $p(n)\in\bbQ[n]$ where $\bbQ$ is the rationals.
It follows from Corollary 1.9.3 in Stanley's text~\cite{sta:ec1} that
$p(n)$ is an integer for all integral $n$ if and only if the
coefficients in the expansion of $p(n)$ using $\cB_0$ are all
integral.  In particular, this is true for $p(n)=p(S;n)$ by our main
theorem.

One might wonder if the coefficients in the $\cB_0$-expansion of
$p(S;n)$ were also nonnegative.  Unfortunately, it is easy to see from
Theorem~\ref{P(m)} that this is not always the case.  However, we
conjecture that $p(S;n)$ can be written as a nonnegative linear
combination of the elements in another basis.  

Throughout this section, let $S$ be a nonempty admissible set of
constants and $m=\max S$.  Let $\CS_{k}$ be the coefficient of
${n-m\choose k}$ in the expansion of $p(S;n)$, so
\[
p(S;n) = \sum_{k=0}^{m-1} \CS_{k} {n-m\choose k}, 
\]
where we know from Theorem~\ref{P(S):th} that $\CS_{m-1}$ is a positive
integer and $\CS_{k}=0$ for $k\geq m$.

\bcon Each coefficient $\CS_{k}$ is a nonnegative integer.  
  \econ

As evidence for this conjecture, we investigate some special cases.
We will first concern ourselves with $\CS_0$.  In the following
results, we use the usual convention that ${a\choose b}=0$ if $b<0$ or
$b>a$.

\ble
\label{bC_0}
For any nonempty set $S$ with constant elements, we have $\CS_0 =0$.
\ele
\prf
If $S$ is not admissible then $\CS_k=0$ for all $k$.  Now suppose $S$ is  admissible and $m=\max S$. 
From our discussion following the proof of Theorem~\ref{P(S):th} about the case $n=i_S=m$, we see that $p(S;m)=0$.
Since our basis is 
${n-m\choose k}, k\ge0,$ we must therefore have 
$
\CS_0 = p(S;m)=0.  
$
\hqedm

Next we consider what happens for peak sets with one element.

\bpr
If $S=\{m\}$ is admissible then
$$
\CS_k =\case{\dil{m-1\choose k}}{if $k\ge1$,}{0}{if $k=0$.\rule{0pt}{15pt}}
$$
\epr
\prf
Using Theorem~\ref{P(m)} and Vandermonde's convolution give
\bea
p(S;n)
&=&-1+{n-1\choose m-1}\\[5pt]
&=&-1+\sum_{k\ge0} {m-1\choose m-k-1}{n-m\choose k}\\[5pt]
&=&\sum_{k\ge1} {m-1\choose k}{n-m\choose k}
\eea
which is what we wished to prove.\hqedm

To deal with peak sets having two elements, we will need the characteristic function $\chi$ which evaluates to $1$ on a true statement and $0$ on a false one.

\bpr
If $S=\{2,m\}$ is admissible then
\beq
\label{2,m:simp}
\CS_k =(m-3){m-2\choose k-1} +(m-2){m-2\choose k}-{m-2\choose k+m-3}.
\eeq
Furthermore, $\CS_k\ge0$ for all $k$.
\epr
\prf
To prove the formula for $\CS_k $, one first shows that
\beq
\label{2,m}
\CS_k=-2{m-2\choose k+m-3}\chi(\mbox{$m$  even})+
\sum_{j=0}^{m-4} (-1)^j  \left[{m-j-2\choose1}-1\right] {m\choose k+j+1}.
\eeq
Since this equality will be generalized in the next proposition, we will provide the details of the proof there.

To simplify this expression, consider the summation part which we will denote by $\Sigma$.  Use the binomial recursion twice on ${m\choose k+j+1}$, each time reindexing the summation to combine terms, to get
\bea
\Sigma
&=&
(m-3){m-1\choose k}+\sum_{j=0}^{m-4} (-1)^j {m-1\choose k+j+1}\\
&=&
(m-3){m-2\choose k-1} +(m-2){m-2\choose k}+(-1)^{m-4}{m-2\choose k+m-3}.
\eea
Adding in the term of~\ree{2,m} containing $\chi$ yields the desired formula.

To show positivity, if $k\ge2$  then the last binomial coefficient in~\ree{2,m:simp} is zero and the result is obvious.  It is also easy to check the case $k=1$, and $k=0$ is Lemma~\ref{bC_0}.  Thus we are done.
\hqedm

We now consider the case of an arbitrary 2-element set.  While we are able to obtain a general summation formula in this case, it does not seem to simplify readily and so we are only able to prove positivity for roughly half the coefficients.

\bth
If $S=\{l,m\}$ is admissible then $\CS_0 =0$, and for $k\ge1$
$$
\CS_k =-2{m-1\choose k+m-l}\chi(\mbox{$m-l$ even})+
\sum_{j=0}^{m-l-2} (-1)^j  \left[{m-j-2\choose l-1}-1 \right]{m\choose k+j+1}.
$$
Furthermore, $\CS_k \ge0$ for $k\ge(m-2)/2$.
\eth
\prf
We first prove the formula for $\CS_k$.  The case $k=0$ is taken care of by Lemma~\ref{bC_0}.  Let $p(n)=p(S;n)$.   For $k\ge1$ we will prove the formula for $\CS_k $ by fixing $l$ and inducting on $m$.

First consider the base case $m=l+2$.  Then $S_1=\{l\}$ and $S_2=\{l,l-1\}$ which is not admissible.  Thus, using Theorem~\ref{P(m)}, Vandermonde's convolution, and the fact that 
$\CS_0 =0$, we see that  equation~\ree{P(S;n)} becomes
\bea
p(n)&=&p_1(m-1){n\choose m-1}-2p_1(n)\\[5pt]
&=&\left[{m-2\choose l-1}-1\right]{n\choose m-1} - 2 \left[{n-1\choose l-1} -1\right]\\[5pt]
&=&\sum_{k\ge1} \left[{m-2\choose l-1}-1\right]{m\choose m-1-k}{n-m\choose k} -2{m-1\choose l-k-1}{n-m\choose k}\\[5pt]
&=&\sum_{k\ge1}\left\{ \left[{m-2\choose l-1}-1\right]{m\choose k+1} - 2{m-1\choose k+m-l} \right\} {n-m\choose k}.
\eea
The coefficient of ${n-m\choose k}$ in this expression agrees with the one given in the theorem when $m=l+2$ and so we done with the base case.

Now consider $m>l+2$.  There are two similar subcases depending on whether $m-l$ is even or odd and so  we will just do the latter. The computations in the base case remain valid except for the fact that $S_2=\{l,m-1\}$ is now admissible and so we need to subtract off the $p_2(n)$ term in equation~\ree{P(S;n)}.  For simplicity, let $a_k$ denote the coefficient of $p_2(n)$ expanded in the basis $\cB_{m-1}$.  Since $m-l-1$ is even we have, by induction,
$$
a_k=-2{m-2\choose k+m-l-1} + \sum_{j=0}^{m-l-3}(-1)^j \left[{m-j-3\choose l-1}-1\right]{m-1\choose k+j+1}
$$
when $k\ge1$ and, as always, $a_0=0$.
To convert to the basis $\cB_m$, we compute
$$
p_2(n)=\sum_{k\ge0}a_k{n-m+1\choose k}
=\sum_{k\ge0}a_k\left[{n-m\choose k-1}+{n-m\choose k}\right]
=\sum_{k\ge0}(a_k+a_{k+1}){n-m\choose k}.
$$
It follows from the previous two displayed equations that the coefficient of ${n-m\choose k}$ in $-p_2(n)$ is
$$
2{m-1\choose k-m-l}-\sum_{j=0}^{m-l-3}(-1)^j \left[{m-j-3\choose l-1}-1\right]{m\choose k+j+2}.
$$
Shifting indices in this last sum and adding in the contribution from the computation for $p(n)$ in the base case completes the induction step. 

The proof of positivity breaks down into two cases depending on the parity of $m-l$.  Since they are similar, we will only present the details when $m-l$ is odd.  It suffices to show that the absolute values of the terms in the sum for $\CS_k $ are weakly decreasing since then each negative term can be canceled into the preceding positive one.    Clearly the term in square brackets is decreasing with $j$.  And because $k+1\ge m/2$ we have that ${m\choose k+j+1}$ is also decreasing by unimodality of the rows of Pascal's triangle.  This completes the proof.
\hqedm

\section{Equidistribution}
\label{e}

 Suppose one considers the distribution of  $\#P(S;n)$ over all possible peak sets $S=\{i_1,\ldots, i_s\}$ with $s$ elements.   We conjecture  that a maximum will occur when the elements of $S$ are as evenly spaced as possible. 

There are two natural probabilistic conjectures which one could make about the peak distribution, assuming a small number of peaks in a long sequence. First, one could guess that the places where peaks occur is an approximation to Poisson process arrivals, and hence locations of the peaks are distributed approximately uniformly over the whole sequence and  are approximately independent. Available evidence points to the alternative conjecture that the peaks have a tendency to repel each other. This phenomenon is found in some random models, e.g. under certain assumptions, eigenvalues of random matrices have a tendency to repel each other. We do not see a direct connection with that model at the technical level, but the repelling nature of peaks invites further exploration.

It will be useful to pass from the set $S$ to the corresponding composition.  A \emph{composition of $n$ into $k$ parts} is a sequence of positive integers $\ka=(\ka_1,\ldots,\ka_k)$ where $\sum_j \ka_j=n$.  We also write 
$\ka=(a^{m_a},b^{m_b},\ldots)$ for the composition which starts with $m_a$ copies of the part $a$, then $m_b$ copies of the part $b$, and so forth.  Given any set $S=\{i_1,\ldots,i_s\}$ of $[n]$ there is a corresponding composition $\ka(S)$ of $n+1$ into $s+1$ parts given by $\ka_j=i_j-i_{j-1}$ for $1\le j\le s+1$ where we let $i_0=0$ and $i_{s+1}=n+1$.  This construction is bijective.  Given any composition $\ka=(\ka_1,\ldots,\ka_{s+1})$ of $n+1$ we can recover $S=\{i_1,\ldots,i_s\}\sbe[n]$ where $i_j=\ka_1+\cdots+\ka_j$ for $1\le j\le s$.   

A composition is \emph{Tur\'an} if $|\ka_a-\ka_b|\le 1$ for all $a,b$.  This terminology is in reference to Tur\'an's famous theorem in graph theory (about maximizing the number of edges in a graph with no complete subgraph of given order) where these compositions play an important r\^ole.   There is another description of Tur\'an compositions which will be useful.  Suppose we wish to form a Tur\'an composition of $n$ with $k$ parts.  Apply the Division Algorithm to write $n=qk+r$ where $0\le r< k$.  Then the desired compositions are exactly those gotten by permuting $k-r$ copies of the part $q$ and $r$ copies of the part $q+1$.   We will call $q$ the \emph{quotient} corresponding to the Tur\'an composition.

\bcon[Equidistribution Conjecture]
First, consider the case when $n,s$ are fixed positive integers.  Then, we conjecture the following hold.
\ben
\item[(a)]  If $S\sbe[n]$ maximizes $\#P(S;n)$ among all subsets with $\#S=s$, then $\ka(S)$ is Tur\'an.
\item[(b)]  The maximizing Tur\'an compositions in (a) are precisely those of the form 
$$
((q+1)^{m_1},q^{m_2},(q+1)^{m_3})
$$
where $q$ is the quotient of $\ka(S)$ and as many of the multiplicities $m_1,m_3$ are positive as possible.  (If there is only one copy of $q+1$, then one of these two multiplicities is zero and the other equals one, and if there are no copies then both multiplicities are zero.)
\een
Next, consider the case when  $n$ is fixed and $s=\#S$ is allowed to vary, then the peak sets maximizing $\#P(S;n)$ over all  $S\sbe[n]$ are the Tur\'an compositions satisfying (b) with the maximum number of $3$'s.
\econ

Note that for $s=1$ this conjecture is true because of Theorem~\ref{P(m)}.  It has also been verified by computer for $n\le 13$.  The part of the conjecture about maximization over all $S\sbe[n]$ is consistent with a result of Kermack and McKendrick~\cite{km:trs} stating that the mean size of a part in all $\ka(S)$ with $S$ admissible is 3.

\section{Peaks and valleys}
\label{pv}

For some applications, it will be useful to know the number of permutations with peaks at $2$ and $n-1$ and a valley at a given position $m$. In the mass redistribution model analyzed in~\cite{bbps:me}, valleys represent the oldest sites, where age is measured since the last mass redistribution. It is a natural question to investigate the relationship between the oldest sites (valleys) and the sites most recently affected by the mass redistribution process (peaks).   

In this section we derive the desired formula.  To set up notation, let
\bea
PV(\pi)&=&\mbox{set of peaks and valleys of permutation $\pi$},\\
PV(i_1,\ldots,i_s;n)&=&\{\pi\in\fS_n\ :\ \mbox{$PV(\pi)=\{i_1,\ldots,i_s\}$ and $i_1$ is a peak}\}.
\eea   	
As usual, we require $i_1<\cdots<i_s$.   Of course, peaks and valleys must alternate.  So $PV(i_1,\ldots,i_s;n)$ also counts permutations  $\pi$ with $PV(\pi)=\{i_1,\ldots,i_s\}$ and $i_1$ being a valley, a fact which will be useful in the sequel.  The definition of admissible is as before.

It is easy to adapt the proof of Theorem~\ref{P(S):th} to this setting, so the demonstration of the next result is omitted.  
\bth
\label{PV(S):th}
If  $S=\{i_1<\cdots<i_s\}$ is admissible then
$$
\#PV(S;n) = q(n) 
$$
where $q(n)$ is a polynomial depending on $S$ such that $q(n)$ is an integer for all integral $n$.  In addition, if  $S$ is a set of constants not depending on $n$ then $\deg q(n) = i_s - 1$ (when $S=\emp$ we have $\deg q(n)=0$). \hqedm
\eth

We now derived a formula for $\#PV(2,m,n-1;n)$ via a sequence of results.  Since the techniques are much like those we have used before, the proofs will only be sketched.

\bpr
If $\{m\}$ is admissible then
$$
\#PV(m;n)={n-1\choose m-1}.
$$
\epr
\prf
We have $\pi\in PV(m;n)$ if and only if $\pi=\pi_L n \pi_R$ where $\pi_L$ is increasing, $\#\pi_L=m-1$, $\pi_R$ is decreasing, and $\#\pi_R=n-m$ and is decreasing.
\hqedm

\bpr
If $\{2,m\}$ is admissible then
$$
\#PV(2,m;n)={n-2\choose m-2} + (m-2){n-1\choose m-1}.
$$
\epr
\prf
If $a_1\ldots a_n\in PV(2,m;n)$ then either $a_1=1$ or $a_m=1$.  In the first case, the number of $\pi$ is given by 
${n-2\choose m-2}$ by the previous proposition.  In the second case, there are ${n-1\choose m-1}$ ways to pick the elements to the left of $1$ and then $m-2$ ways to pick $a_1$.
\hqedm

\bpr
If $\{2,m,n-1\}$ is admissible then
$$
\#PV(2,m,n-1;n)=2(n-1)\left[{n-4\choose m-2}+(m-2){n-3\choose m-1}\right].
$$
\epr
\prf
By symmetry, it suffices to count the number of $a_1\ldots a_n\in PV(2,m,n-1;n)$ where $a_{n-1}=n$ and double.  There are $n-1$ ways to choose $a_n$.  And using the previous proposition, we see  that the number of ways to pick the remaining elements is given by the expression in the square brackets.
\hqedm

\section{Fixing the number of peaks and future research}
\label{fnp}

We now use our theorems to prove a result already in the
literature.  In general, there does not seem to be a simple explicit
formula for the number $f(s,n)$ of permutations in $\fS_n$ with $s$
peaks, 
see sequence A008303 in the Online Encyclopedia of
Integer Sequences (OEIS). However, David and Barton~\cite[p.\ 163]{db:cc}
give the recurrence
\[
f(s,n) = (2s +2) f(s,n-1) + (n-2s)f(s-1,n-1)
\]
with the initial conditions that $f(0,n)=2^{n-1}$ and $f(s,n)=0$
whenever $s\geq \frac{n}{2}$.  In addition, 
for small $s$, one can write
down an explicit expression for $f(s,n)$.  In fact, the sequence
$f(1,n)$ appears 
as sequence A000431 in the OEIS where the following result is attributed to Mitchell
Harris.  \bpr For $n\ge1$
$$
f(1,n)=2^{2n-3}-n2^{n-2}.
$$
\epr
\prf
Using Theorems~\ref{P(S):th} and~\ref{P(m)} we have
$$
f(1,n)=\sum_{m=2}^{n-1} \left[{n-1\choose m-1}-1\right] 2^{n-2}
=2^{n-2}\sum_{m=1}^{n} \left[{n-1\choose m-1}-1\right] =2^{n-2}(2^{n-1}-n)
$$
which multiplies out to the formula we want.
\hqedm

We suggest the following problem for the interested reader.  In the current work we have only considered peak sets containing
certain linear functions of $n$.  It would be interesting for see what happens if $S$ contains other functions of $n$, for example, $\sqrt{n}$ or $\ln{n}$.

\medskip

{\it Acknowledgments.}  We would like to thank Andrew Crites and
Soumik Pal for stimulating discussions about peak sets.  We would also
like to thank Andrew Sills for information about hypergeometric
series.  We credit the OEIS for assisting us with our research on this project.


\begin{thebibliography}{10}

\bibitem{ano:nrp}
Marcelo Aguiar, Kathryn Nyman, and Rosa Orellana.
\newblock New results on the peak algebra.
\newblock {\em J. Algebraic Combin.}, 23(2):149--188, 2006.

\bibitem{BHvW}
Louis~J. Billera, Samuel~K. Hsiao, and Stephanie van Willigenburg.
\newblock Peak quasisymmetric functions and {E}ulerian enumeration.
\newblock {\em Adv. Math.}, 176(2):248--276, 2003.

\bibitem{bbps:me}
Sara Billey, Krzysztof Burdzy, Soumik Pal, and Bruce~E. Sagan.
\newblock On meteors and earthworms.
\newblock Forthcoming (2012).

\bibitem{bcmyy:vpp}
Pierre Bouchard, Hungyung Chang, Jun Ma, Jean Yeh, and Yeong-Nan Yeh.
\newblock Value-peaks of permutations.
\newblock {\em Electron. J. Combin.}, 17(1):Research Paper 46, 20, 2010.

\bibitem{db:cc}
F.~N. David and D.~E. Barton.
\newblock {\em Combinatorial chance}.
\newblock Hafner Publishing Co., New York, 1962.

\bibitem{fv:psl}
Jean Fran{\c{c}}on and G{\'e}rard Viennot.
\newblock Permutations selon leurs pics, creux, doubles mont\'ees et double
  descentes, nombres d'{E}uler et nombres de {G}enocchi.
\newblock {\em Discrete Math.}, 28(1):21--35, 1979.

\bibitem{GKP}
Ronald~L. Graham, Donald~E. Knuth, and Oren Patashnik.
\newblock {\em Concrete mathematics}.
\newblock Addison-Wesley Publishing Company, second edition, 1994.
\newblock A foundation for computer science.

\bibitem{km:trs}
W.~O. Kermack and A.~G. McKendrick.
\newblock Tests for randomness in a series of numerical observations.
\newblock {\em Proc. Roy. Soc. Edinburgh}, 57:228--240, 1937.

\bibitem{ma:dpe}
Shi-Mei Ma.
\newblock Derivative polynomials and enumeration of permutations by number of
  interior and left peaks.
\newblock {\em Discrete Math.}, 312(2):405--412, 2012.

\bibitem{nym:pas}
Kathryn~L. Nyman.
\newblock The peak algebra of the symmetric group.
\newblock {\em J. Algebraic Combin.}, 17(3):309--322, 2003.

\bibitem{pet:epp}
T.~Kyle Petersen.
\newblock Enriched {$P$}-partitions and peak algebras.
\newblock {\em Adv. Math.}, 209(2):561--610, 2007.

\bibitem{schocker}
Manfred Schocker.
\newblock The peak algebra of the symmetric group revisited.
\newblock {\em Adv. Math.}, 192(2):259--309, 2005.

\bibitem{sta:ec1}
Richard~P. Stanley.
\newblock {\em Enumerative combinatorics. {V}olume 1}, volume~49 of {\em
  Cambridge Studies in Advanced Mathematics}.
\newblock Cambridge University Press, Cambridge, second edition, 2012.

\bibitem{ste:epp}
John~R. Stembridge.
\newblock Enriched {$P$}-partitions.
\newblock {\em Trans. Amer. Math. Soc.}, 349(2):763--788, 1997.

\bibitem{str:eap}
Volker Strehl.
\newblock Enumeration of alternating permutations according to peak sets.
\newblock {\em J. Combinatorial Theory Ser. A}, 24(2):238--240, 1978.

\bibitem{ws:pen}
Di~Warren and E.~Seneta.
\newblock Peaks and {E}ulerian numbers in a random sequence.
\newblock {\em J. Appl. Probab.}, 33(1):101--114, 1996.

\end{thebibliography}
\end{document}